\documentclass[a4paper,10pt]{article}
\usepackage{stmaryrd}
\usepackage{amsfonts}
\usepackage{bbm}
\usepackage{amscd}
\usepackage{mathrsfs}
\usepackage{latexsym,amssymb,amsmath,amscd,amscd,amsthm,amsxtra,xypic}
\usepackage[dvips]{graphicx}

\newtheorem{thm}{Theorem}[section]
\newtheorem{lem}[thm]{Lemma}
\newtheorem{cor}[thm]{Corollary}
\newtheorem{pro}[thm]{Proposition}

\newtheorem{defi}[thm]{Definition}

\newcommand{\lon }{\,\rightarrow\,}
\newcommand{\be }{\begin{eqnarray*}}
\newcommand{\ee }{\end{eqnarray*}}

\newcommand{\ga}{\Gamma ({\mathcal A})}

\newcommand{\cm}{C^{\infty}(M)}

\newcommand{\rr}{\Lambda^{\#}}

\newcommand{\rt}{\Lambda}

\newcommand{\f}[1]{{\mathfrak{#1}}}

\newcommand{\defbe}{\stackrel{\triangle}{=}}
\newcommand{\pf}{\noindent{\bf \quad\quad Proof.}\ }

\newcommand{\inverse}{^{-1}}

\newcommand{\Real}{{{\mathbb R}}}

\newcommand{\huaA}{{\mathcal{A}}}

\newcommand{\huaD}{{\mathcal{D}}}

\newcommand{\CWM}{C^{\infty}(M)}
\newcommand{\CWGPV}{C^{\infty,G}(P,V)}

\newcommand{\adjoint}{{ad}}

\newcommand{\XM}{{\mathcal{X}}(M)}

\newcommand{\set}[1]{\left\{#1\right\}}

\newcommand{\frkg}{\mathfrak{g}}

\newcommand{\frkx}{\mathfrak{x}}

\newcommand{\frkD}{\mathfrak{ D}}

\newcommand{\frkL}{\mathfrak{ L}}

\def\gpd{\,\lower1pt\hbox{$\longrightarrow$}\hskip-.24in\raise2pt
         \hbox{$\longrightarrow$}\,}

\def\qed{\hfill ~\vrule height6pt width6pt depth0pt}

\newcommand{\Rep}{{\frkL}}

\newcommand{\LieG}{{\frkg}}
\newcommand{\lieG}{Lie(G)}

\newcommand{\CDO}{{\mathcal{CDO}}}
\newcommand{\DiffOper}{{\frkD}}
\begin{document}
\title{
{The Cohomology of Transitive Lie Algebroids
\thanks
 {
 Partially supported by NSFC(19925105) and
CPSF(20060400017), this article appeared in \emph{NONCOMMUTATIVE
GEOMETRY AND PHYSICS
 2005,
  Proceedings of the International Workshop Sendai-Beijing Joint Workshop}
  (Page 109-127, World Scientific Publishing Company,  Singapore).
 }
} }
\author{ Z. Chen  and ~ Z.-J. Liu \\
Department of Mathematics and LMAM\\ Peking University,
Beijing 100871, China\\
          {\sf email: chenzhuo@math.pku.edu.cn,\quad liuzj@pku.edu.cn } }

\date{}
\footnotetext{{\it{Keywords}}: Lie algebroid, cohomology, Lie
bialgebroid.}

\footnotetext{{\it{MSC}}:  Primary 17B65. Secondary 18B40, 58H05. }

\maketitle  

\begin{abstract}
For a transitive Lie algebroid $\huaA$ on a connected manifold $M$
and its a representation on a vector bundle $F$, we study the
localization map $\Upsilon^1$: $H^1(\huaA,F)\lon H^1(L_x,F_x)$,
where $L_x$ is the adjoint algebra at $x\in M$. The main result in
this paper is that: $Ker\Upsilon^1_x=Ker(p^{1*})=H^1_{deR}(M,F_0)$.
Here $p^{1*}$ is the lift of $H^1(\huaA,F)$ to its counterpart over
the universal covering space $\widetilde{M}\stackrel{p}{\lon} M$ and
$H^1_{deR}(M,F_0)$ is the $F_0=H^0(L,F)$-coefficient deRham
cohomology. We apply these results to study the associated vector
bundles   to   principal fiber bundles and the structure of
transitive Lie bialgebroids.
\end{abstract}
\section{Introduction}

The theory of Lie algebroids is one of  important fields in modern
differential geometry, which gives an unified way to study Lie
algebras and the tangent bundle of a manifold. Its global version is
Lie groupoid. Please see \cite{first} for a detailed introduction to
the theory of  Lie algebroids and Lie groupoids as well as
\cite{WeinGroupoids} by Weinstein for their applications in Poisson
geometry. Moreover, in \cite{Connes} Connes also points out that
groupoids play an essential role in non-commutative geometry.

The purpose of this paper is to study the cohomology of Lie
algebroids, which is a basic topic in this field and has wide
applications in physics and other fields of mathematics (e.g., see
\cite{Evens}, \cite{Crabowski} and \cite{Kapusting}), as well as its
applications for  associated vector bundles with respect to some
principal fiber bundle and the structure of transitive  Lie
bialgebroids. The notion of a Lie bialgebroid was introduced by
Mackenzie and Xu in \cite{MackenzieX:1994} as a natural
generalization of that of a Lie bialgebra, as well as the
infinitesimal version of  Poisson groupoids introduced by Weinstein
\cite{Weinstein:1988}. It has been shown that much of the theory of
Poisson groups and Lie bialgebras can be similarly carried out  in
this general context. It is therefore a basic task to  study the
structure of Lie bialgebroids. In particular,  it is very
interesting to figure out what special features a Lie bialgebroid,
in which the Lie algebroid structure is transitive, would have.

The paper is organized as follows. In Section 2, first we give a
brief introduction to basic notion and then,
 for a transitive Lie algebroid $\huaA$ and its representation
 on a vector bundle $F$, we define  a morphism of cohomology
 groups, called the localization map, and prove that the
 Lie algebroid $1$-cohomology  is totally determined by
the $1$-cohomology of its adjoint Lie algebra under some topological
condition.  In Section 3 we study the kernel of the localization map
mentioned above. To do it, the connection and parallelism  are
used(see \cite{CL} for more details).
 As an equivalent statement of the main theorem, we describe some
conditions to make a Lie algebroid $1$-cocycle to be coboundary. In
section 4, we study some properties of the cohomology of associated
vector bundles with respect to some principal fiber bundle.
 In Section 5,   we recall some results on the  structure of transitive Lie bialgebroids in
\cite{CL} as an another  application of our localization theory.

\section{Preliminaries}
Throughout the paper we suppose that any smooth manifold $M$ under
consideration is connected.


\subsection{Lie algebroids and its representations} First we introduce some  basic concept used below:
\begin{defi}\label{Def:LieAlgebroid}
A Lie algebroid ${\huaA}$ with base space $M$, is a (real) vector
bundle over $M$, together with a bundle map $\rho: {\huaA}\lon TM$,
called the anchor, and there is a (real) Lie algebra structure
$[~\cdot~,~\cdot~]_{{\huaA}}$ on $\Gamma({\huaA})$ satisfying the
following conditions:
\begin{itemize}
\item[1.] The induced map $\rho: \Gamma({\huaA})\lon \XM$ is a Lie algebra morphism;
\item[2.] for all $f\in \CWM$, $A,B\in \Gamma({\huaA})$, the Leibnitz law holds. That is,
$$[A,fB]_{{\huaA}}=f[A,B]_{{\huaA}}+(\rho(A)f)B.$$
\end{itemize}
\end{defi}
In this paper, by $({\huaA}$, $[~\cdot~,~\cdot~]_{{\huaA}},\rho)$ we
denote a Lie algebroid ${\huaA}$ with anchor map $\rho$ :
${\huaA}\lon TM$ and Lie bracket $[~\cdot~,~\cdot~]_{{\huaA}}$ on
$\Gamma({\huaA})$. For a transitive Lie algebroid
$({\huaA},[~\cdot~,~\cdot~]_{{\huaA}},\rho)$ over $M$ (i.e., the
anchor $\rho$ is surjective), the \emph{Artiyah sequence} is as
follows:
\begin{equation}\label{SequenceAtiyah}
 0\lon  L\stackrel{i}{\longrightarrow} {\huaA}
\stackrel{\rho}{\longrightarrow} TM\lon  0,
\end{equation}
where $L=Ker(\rho)$ is called  the {\em adjoint bundle} of
${\huaA}$. In fact, $L$ is a Lie algebra bundle. The fiber type can
be taken as the Lie algebra $\f{g}=L_p$ at any point $p\in M$. By
$[~\cdot~,~\cdot~]_L$, we denote the fiber-wise bracket on $L$ (see
 \cite{first}).

Let $F$ be a vector bundle over $M$. By $\CDO(F)$
 we denote the \emph{bundle of covariant differential operators} of $F$.
 (see III in  \cite{first} and refs. \cite{Mackenzie:1995,KSK}
 for more details. In the last text, the authors use ${\huaD}(F)$
instead of $\mathcal{CDO}(F)$). Each element $D\in$ $\CDO(F)_{x}$ is
an operator $D_x: \Gamma(F)\lon F_x$, corresponding to a unique
$X\defbe \sigma(D_x)\in T_xM$, such that
\begin{equation}
D_x(f\mu)=X(f)\mu(x)+f(x)D_x\mu,\quad\forall f\in
C^{\infty}(M),\mu\in\Gamma(F).
\end{equation}
Actually, $(\CDO(F),[~\cdot~,~\cdot~]_{\CDO(F)},\sigma)$ is a
transitive Lie algebroid over $M$, with the bracket of commutator
$[~\cdot~,~\cdot~]_{\CDO(F)}$. That is, for two
$D_1,D_2\in\Gamma(\CDO(F))$,
$$[D_1,D_2]_{\CDO(F)}\defbe D_1\circ D_2-D_2\circ D_1,$$
is also an element of $\Gamma(\CDO(F))$. Moreover, the corresponding
vector field
$$\sigma[D_1,D_2]_{\CDO(F)}=[\sigma(D_1),\sigma(D_2)].$$
Obviously, the adjoint bundle of this Lie algebroid is
$End(F)$=$Ker(\sigma)$.

\begin{defi}
A representation of a Lie algebroid
$({\huaA},[~\cdot~,~\cdot~]_{{\huaA}},\rho)$ on a vector bundle
$F\lon M$ is a vector bundle map $\Rep: {\huaA}\lon \CDO(F)$, which
is also a Lie algebroid morphism. That is,
\begin{itemize}
\item[1)] $\forall A\in {\huaA}_x$, $x\in M$, $\rho(A)=\sigma(\Rep(A))$;
\item[2)] $\forall A,B\in\Gamma({\huaA})$,
$\Rep[A,B]_{{\huaA}}=[\Rep(A),\Rep(B)]=\Rep(A)\circ\Rep(B)-
\Rep(B)\circ\Rep(A)$.
\end{itemize}
\end{defi}
We also denote $\Rep(A)$ by $\Rep_A$. For example, when ${\huaA}$ is
transitive, one can define the adjoint representation of ${\huaA}$
on $L$: for each $A\in{\huaA}_x$, define
\begin{equation}\label{adjointRep}
\adjoint_A: \Gamma{(L)} \lon L_x, \quad \mu\mapsto
[A,\mu]_{{\huaA}}(x), \quad\forall \mu\in \Gamma({L}).
\end{equation}
 Here, we need
to extend $A$ to be a local section $\tilde{A}$ of ${\huaA}$ near
$x$, and then the value of $[A,\mu]_{{\huaA}}(x)$ is defined to be
$[\tilde{A},\mu]_{{\huaA}}(x)$. Note that, it does not depend on the
choice of the extension of $A$ near $x$. Usually $\adjoint_A$ is
written as $[A,~\cdot~]_{{\huaA}}$.

As usual, a representation $\Rep$ of ${\huaA}$ on $F$ defines  the
Chevalley complex $(C^k({\huaA},F),$ ${\DiffOper})$ and cohomology
groups $H^k({\huaA},{F})$ $=$ $Ker({\DiffOper})/Im({\DiffOper})$
(see  \cite{Crainic}). In detail, we call a series of $\CWM$-modules
$$
C^k({\huaA},{F})\defbe \set{\mbox{bundle maps }~~ \Omega:
\wedge^k{\huaA}\lon {F} }, \ (k\geqslant 1), $$ and
$C^0({\huaA},{F})=\Gamma({F})$ the cochain space of $\Rep$. The
coboundary operators (also called the differentials)
$\DiffOper=\DiffOper^k: C^k({\huaA},{F})\lon C^{k+1}({\huaA},{F})$
are defined in the traditional way and satisfy $\DiffOper^2=0$. For
this complex $(C^k({\huaA},V), \DiffOper)$, the corresponding
cohomology groups are
$$H^k({\huaA},{F})={Ker}(\DiffOper^k)/Im(\DiffOper^{k-1}).\quad (k\geqslant  0)$$
We adopt the convention that $H^0({\huaA},{F})=Ker(\DiffOper^0)$. In
particular,  a closed 0-cochain is a smooth section of $\Gamma(F)$,
say $\nu$, satisfying $\Rep_A \nu=0$, $\forall A\in {\huaA}$. The
group $H^0({\huaA},{F})$ is the collection of all such closed
0-cochains. A 1-cochain $\Omega\in C^1({\huaA}, {F})$ is a bundle
map from ${\huaA}$ to $F$ . It is called closed, or a
1-\emph{cocycle}, denoted by $\Omega\in D^1({\huaA},F)$, if
\begin{equation}\label{EqtCocycleOmega}
\Omega[A,B]_{{\huaA}} =\Rep_{A}(\Omega(B))-\Rep_{B}(\Omega(A)),
\quad\forall A,B\in\Gamma({\huaA}).
\end{equation}
Especially, we call $\Omega$ a coboundary, denoted by $\Omega\in
B^1({\huaA},F)$, if $\Omega={\DiffOper}\mu$, i.e.,
$$
\Omega(A)=\Rep_{A}(\mu),\quad\forall A\in\Gamma({\huaA})
$$
for some $\mu\in\Gamma(F)$. In what follows,  cochains are simply
called chains. Two 1-cocycles are called homologic, if their
substraction is a coboundary. The group $H^1({\huaA},{F})$ are the
quotient group of all 1-cocycles in sense of homological
equivalence:
$$
H^1({\huaA},{F})\defbe D^1({\huaA},F)/B^1({\huaA},F).
$$
We usually write the equivalence class of a 1-cocycle $\Omega$ by
$[\Omega]$.

If the Lie algebroid degenerates to a Lie algebra, then the above
construction of cohomology groups returns to that of the Lie
algebras.


\subsection{The localization of 1-cohomology of transitive Lie
algebroids} We choose an arbitrary point $x\in M$. For the Lie
algebra $(L_x,[~\cdot~,~\cdot~]_{L_x})$, $\Rep$ induces a
representation $\underline{\Rep}$ of $L_x$ on $F_x$. In fact, for
any $u\in L_x$, we define $ \underline{\Rep}_{u}(\mu)\defbe
\Rep_{u}(\overline{\mu}) $, here $\mu\in F_x$,
$\overline{\mu}\in\Gamma(F)$ is a locally smooth extension of $\mu$.
This is well defined, because $\Rep_{u}(f\mu)=f(x)\Rep_{u}(\mu)$,
$\forall f\in\CWM$. Since $\Rep$ is a Lie algebra morphism, so is
$\underline{\Rep}$: $L_x\lon End(F_x)$ and hence $\underline{\Rep}$
is indeed a representation.

Consider the Chevalley complex $C^\bullet(L_x,F_x)$, where $F_x$ is
regarded as an $L_x$-module via the representation $\underline\Rep$.
We denote the set of closed chains in $C^k(L_x,F_x)$ by
$D^k(L_x,F_x)$, and the corresponding Chevalley cohomology groups by
$H^k(L_x,F_x)$.

We need the following theory of localizations and some results in
\cite{partone}.
\begin{pro}\label{Thm:Upsilon}
Let $\Rep: {\huaA}\lon \CDO(F)$ be a representation of a Lie
algebroid $({\huaA},[ \cdot , \cdot ]_{{\huaA}},\rho)$ on $F\lon M$.
For any $x\in M$, there is a natural morphism $\Upsilon_x^k:$
$H^k({\huaA},{F})\lon H^k(L_x,F_x)$ defined by
$$\Upsilon_x^k([\Omega])\defbe [\underline{\Omega}],
\quad\forall \Omega\in C^k({\huaA}, {F}).$$ Here
$\underline{\Omega}$ is the limitation of $\Omega$ on $\wedge^k L_x
\subset \wedge^k{\huaA}_x$.
\end{pro}
We call the group morphism $\Upsilon^k_x$ defined above the {\bf
localization} of the homology group $H^k({\huaA},F)$ at $x$. For
this $\Upsilon$, we have the following claims.

\begin{thm}\label{Thm:UpsilonEqual}
Given a representation of a {\bf transitive} Lie algebroid
$({\huaA},[~\cdot~,~\cdot~]_{{\huaA}},\rho)$
 on $F\lon M$, $\Rep: {\huaA}\lon \CDO(F)$, then we have
\begin{itemize}
\item[1)] For any $x$, $y\in M$, there exits an isomorphism $J$:
$H^1(L_y,{F_y})\lon H^1(L_x,{F_x})$, such that the following diagram
commutes
\begin{equation}\label{Eqt:UpsilonCommute0}
\xymatrix@R=0.5cm{
                &         H^1(L_y,{F_y}) \ar[dd]^{J}     \\
H^1({\huaA},{F})\quad\quad
 \ar[ur]^{{\Upsilon_y^1}} \ar[dr]_{{\Upsilon_x^1}}                 \\                &         H^1(L_x,{F_x})                 }
\end{equation}
\item[2)]
If $M$ is simply connected, or $H^0(L_x,{F_x})=0$, then the
localization $\Upsilon^1_x$: $H^1({\huaA},{F})$ $\lon$
$H^1(L_x,{F_x})$ is an injection.
\end{itemize}
\end{thm}
In general, the isomorphism $J$ in 1) depends on the choice of a
path from $x$ to $y$ in $M$. It is shown in  \cite{partone} that,
although the isomorphism $J$: $H^1(L_y,{F_y})\cong H^1(L_x,{F_x})$
is not naturally defined, the two subgroups $Im\Upsilon_y^1$ and
$Im\Upsilon_x^1$ are naturally isomorphic, under the condition of 2)
of this theorem. In this paper, we are going to prove in Corollary
\ref{CorImYCong} that even without this condition, the conclusion
also holds.

\section{Kernel of the localization map $\Upsilon^1$}
An interesting problem is that, if $M$ is not simply connected, and
$H^0(L_x,{F_x})\neq 0$, then what is the kernel of $\Upsilon_x^1$?
In this section, we give the answer to this problem in two different
ways. Please see the following two
equations 
(\ref{Eqt:KerUpsilon1}) and (\ref{Eqt:KerUpsilon2}).
\begin{lem}
Let $p: \widetilde{M}\lon M$ be a covering map. The for any vector
bundle $F\lon M$, we have
$$\CDO(p^!{F})\cong p^!\CDO(F).$$
\end{lem}
\pf For each pair $(z,D_0)\in p^!\CDO(F)_z$, where $z\in
\widetilde{M}$,  $x=p(z)$, $D_0\in\CDO(F)_x$, we define
$\phi(z,D_0)\in \CDO(p^!F)_z$ as follows: given an arbitrary
$\lambda\in\Gamma(p^!F)$, we find a decomposition near $z$:
$$\lambda=\sum_i f_i\mu_i,\quad
\mbox{where }\ f_i\in C^\infty{(\widetilde{M})},\mu_i\in
\Gamma(F),$$ and we then set
$$
\phi(z,D_0)\lambda\defbe \sum_i (Z(f_i)\mu_i(x)+f_i(z)D_0(\mu_i)),
$$
where $Z\in T_z\widetilde{M}$ is the unique tangent vector
satisfying $p_*(Z)=\sigma(D_0)$. It is easy to prove that this
definition does not depend on the choice of decompositions of
$\lambda$. In this way, we obtain $\phi: p^!\CDO(F)\lon \CDO(p^!F)$,
which is obviously an injection.

Conversely, for any $z\in \widetilde{M}$, $x=p(z)$, and $D\in
\CDO(p^!F)_z$, $\sigma(D)=Z\in T_z\widetilde{M}$, $D$ induces
$D_0\in\CDO(F)_x$, such that $\sigma(D_0)=p_*(Z)$, $\phi(z,D_0)=D$.
In fact, each section $\mu\in \Gamma(F)$ can be naturally regarded
as $\mu\in\Gamma(p^!F)$. Let $D(\mu)=\nu\in (p^!F)_z=F_x$ and we
define $D_0(\mu)=\nu$. Therefore, for the preceding $\lambda
\in\Gamma(p^!F)$, it is not hard to see
$$
D(\lambda)=\sum_i (Z(f_i)\mu_i(x)+f_i(z)D(\mu_i))
=\phi(z,D_0)(\lambda).
$$
This shows that $\phi$ is also surjective. So $\phi$ is indeed an
isomorphism from $p^!\CDO(F)$ to $\CDO(p^!F)$ . \qed

\begin{lem}\label{Lem:PullBackHuaA}
Suppose that a Lie algebroid $(\huaA,[~\cdot~,~\cdot~]_\huaA,\rho)$
has a representation on $F\lon M$, $\Rep: \huaA\lon \CDO(F)$. Let
$p: \widetilde{M}\lon M$ be a covering map. Then the pull back
bundle $\widetilde{\huaA}=p^!\huaA$ is also a Lie algebroid over
$\widetilde{M}$, and it has an induced representation on
$\widetilde{F}=p^!F$, $\widetilde{\Rep}:$ $\widetilde{\huaA}\lon
\CDO(\widetilde{F})$, such that the following diagram commutes.
\begin{equation}
\begin{CD}
\widetilde{\huaA} @> {\widetilde{\Rep}}>>
 ~~~~\ \ \CDO(\tilde{F})
\cong p^!\CDO(F)\\
{p^!}@VV V {p^!} @VVV\\
\huaA @> {\Rep}>>   ~~~~\ \ \CDO(F).
\end{CD}
\end{equation}
In other words, for $z\in \widetilde{M}$, $x=p(z)$, $A\in \huaA_x$,
one has
$$
\widetilde{\Rep}{(x,A)}=(x,\Rep(A)).
$$
Moreover, if $\huaA$ is transitive, then so is $\widetilde{\huaA}$.
\end{lem}
We omit the proofs of this lemma and the following two theorems.

\begin{thm}\label{Thm:KerUpsilon11}
With the same assumptions as in Lemma 
\ref{Lem:PullBackHuaA}, we have a group morphism
$$p^{k*}: H^k(\huaA,F)\lon H^k(\widetilde{\huaA},\widetilde{F}).$$
\end{thm}

\begin{thm}\label{Thm:KerUpsilon12}
With the same assumptions as in Lemma 
\ref{Lem:PullBackHuaA}, we have a commute diagram
\begin{equation}\left.
\begin{array}{ccc}
H^k(\huaA,{F}) & \stackrel{p^{k*}}{\longrightarrow} &
H^k(\widetilde{\huaA},\widetilde{F}) \\
{\scriptstyle {\Upsilon_x^k}}{\downarrow} & ~&
\downarrow {\scriptstyle \widetilde{\Upsilon}_{\widetilde{x}}^k}\\
H^k(L_x,{F_x})&=& H^k(L_{\widetilde{x}},{F_{\widetilde{x}}}).
\end{array}\right.
\end{equation}
Here $\widetilde{\Upsilon}_{\widetilde{x}}^k$ is the localization of
$H^k(\widetilde{\huaA},\widetilde{F})$ at some $\widetilde{x}\in
\widetilde{M}$ with $p(\widetilde{x})=x$. And one naturally regards
$H^k(L_x,F_x)=H^k(L_{\widetilde{x}},{F_{\widetilde{x}}})$.
\end{thm}

\begin{cor}\label{Cor:KerUpsilon1}
If $\huaA$ is a transitive Lie algebroid and $p$: $\widetilde{M}\lon
M$ is a universal covering, then
\begin{equation}\label{Eqt:KerUpsilon1}
Ker(\Upsilon_x^1)=Ker(p^{1*}).
\end{equation}
\end{cor}
\pf By Theorem 
\ref{Thm:KerUpsilon12},
$\Upsilon_x^1=\widetilde{\Upsilon}_{\widetilde{x}}^1\circ p^{1*}$.
By 2) of Theorem \ref{Thm:UpsilonEqual}
and $\widetilde{M}$ being simply connected, we know
$\widetilde{\Upsilon}_{\widetilde{x}}^1$ is an injection. So we get
(\ref{Eqt:KerUpsilon1}). \qed

The conclusion of Equation 
(\ref{Eqt:KerUpsilon1}) of course describes  $Ker(\Upsilon_x^1)$,
but it has no relationship with the group $H^0(L_x,F_x)$. We now
give another description of $Ker(\Upsilon_x^1)$. For the Lie
algebroid $(\huaA,[~\cdot~,~\cdot~]_\huaA,\rho)$ and its
representation on $F\lon M$, $\Rep: \huaA\lon \CDO(F)$, and any
$x\in M$, we consider a sub vector space
$$F_{0x}=\set{\nu\in F_x| \Rep_u(\nu)=0,\forall u\in L_x}
=H^0(L_x,F_x).$$
Then by Theorem \ref{Thm:UpsilonEqual}, 
when $\huaA$ is transitive,
$$F_0=H^0(L,F)\subset F$$
is a sub vector bundle. Since for each $u\in \Gamma(L)$,
$A\in\Gamma(\huaA)$, $\nu\in\Gamma(F_0)$, we have
$$
\Rep_u(\Rep_A\nu)=\Rep_{[u,A]_\huaA}\nu-\Rep_A(\Rep_u\nu)=0,
$$
and hence $\Rep_A\nu\in\Gamma(F_0)$. So we have an induced
representation of $\huaA$ on $F_0$, also denoted by $\Rep$.

Meanwhile, $\Rep$ induces a representation of $TM$ on $F_0$, denoted
by $\bar{\Rep}: TM\lon \CDO(F_0)$. In fact, for each $X\in T_xM$,
and an arbitrary $A\in\huaA_x$, $\rho(A)=X$, we set
$$\bar{\Rep}_X\nu\defbe \Rep_A \nu, \quad\forall \nu\in \Gamma(F_0).$$
Obviously this definition does not depend on the choice of $A$, and
$\bar{\Rep}$ is well defined. A representation of the tangent bundle
$TM$ on $F_0$ is also referred as a \textit{flat connection} of
$F_0$. We will call $\bar{\Rep}$ the \textbf{reduced (flat)
connection} of $F_0$ coming from the represention $\Rep$. Now,
elements of
$$C^k(TM,F_0)=Hom(\wedge^k(TM),F_0) \quad( \mbox{ with } C^0(TM
,F_0)=\Gamma(F_0))$$ are also called the $F_0$-coefficient
$k$-forms. With the usual exterior differential operator $d:
C^k(TM,F_0)\lon C^{k+1}(TM,F_0)$, $C^{\bullet}(TM,F_0)$ is known as
the $F_0$-coefficient de Rham complex. Especially, $D^1(TM,F_0)$ is
the kernel of $d: C^1(TM,F_0)\lon C^{2}(TM,F_0)$ and
$$
H^1_{deR}(M,F_0)=D^1(TM,F_0)/ d(\Gamma(F_0))
$$
is the first $F_0$-coefficient \textbf{de Rham cohomology} of $M$.
We are going to prove that this group is just the kernel of
$\Upsilon^1_x$
(Theorem \ref{Thm:KerUpsilonSecond}).

\begin{lem}\label{Lem:Correspondence}
There is a one-one correspondence between
$$D^1(\huaA,F_0)_{0}=
\set{\Omega\in D^1(\huaA,F_0)| \Omega(v)=0,\forall v\in L}$$ and
$D^1(TM,F_0)$, denoted by $\Omega\mapsto \bar{\Omega}$. For the
inversion map, we denote
$${\theta}\in D^1(TM,F_0)\ \mapsto\ {\rho^*\theta}\in D^1(\huaA,F_0)_{0}.$$
More over, this map induces an injection of $H^1_{deR}(M,F_0)$ into
$H^1(\huaA,F)$: $[\theta]\mapsto [{\rho^*\theta}]$.
\end{lem}
\pf Given any $\Omega\in D^1(\huaA,F_0)_0$, we define $\bar{\Omega}$
to be a map sending $X\in T_xM$ to $\bar{\Omega}(X)\defbe
\Omega(A)$, where $A\in \huaA_x$, $\rho(A)=X$. This is of course
well defined. It is also easy to check
\begin{equation}\label{Eqt:temp1}
\bar{\Rep}_X\bar{\Omega}(Y)-\bar{\Rep}_Y\bar{\Omega}(X)
-\bar{\Omega}([X,Y])=0, \quad\forall X,Y\in \XM.
\end{equation}
I.e., $\bar{\Omega}\in D^1(TM,F_0)$.

On the other hand, given any ${\theta}\in D^1(TM,F_0)$
satisfying Equation 
(\ref{Eqt:temp1}), we can define
$${\rho^*\theta}(A)={\theta}(\rho(A)),\quad \forall A\in \huaA.$$
${\rho^*\theta}$ naturally satisfies ${\rho^*\theta}|_{L}=0$, and it
is also a cocycle. It is just by the definitions to see that the map
$${\rho^*(\cdot)}: H^1_{deR}(M,F_0)
\lon H^1(\huaA,F_0); \quad [\theta]\mapsto [{\rho^*\theta}],$$ is an
injection of cohomology groups. But as we shall see in the following
exact
Sequence 
(\ref{LongSequenceA}) that $H^1(\huaA,F_0)$ is embedded into
$H^1(\huaA,F)$. So we conclude that $H^1_{deR}(M,F_0)$ can be
embedded into $H^1(\huaA,F)$ via $[\rho^*(\cdot)]$. \qed

Let $\overline{F}=F/F_0$ be the quotient bundle. Now one obtains an
exact sequence
\begin{equation}\label{SequenceF}
 0\lon  F_0\stackrel{i}{\longrightarrow} F
\stackrel{j}{\longrightarrow} \overline{F}\lon  0.
\end{equation}
The algebroid $\huaA$ has an induced representation on
$\overline{F}$, in an obvious sense and also denoted by $\Rep$:
$$
\Rep_A [\mu]\defbe [\Rep_A \mu],\quad\forall \mu\in \Gamma(F).
$$

Now, we get an exact sequence of complexes
\begin{equation}\label{SequenceCAF}
 0\lon C^k(\huaA, F_0)\stackrel{i}{\longrightarrow} C^k(\huaA,F)
\stackrel{j}{\longrightarrow} C^k(\huaA,\overline{F})\lon  0.
\end{equation}
Here $i$, $j$ are both cochain maps. So we have the following long
exact sequence (the Mayer-Vietoris) of cohomology groups
\begin{eqnarray}\nonumber
&& 0\longrightarrow H^0(\huaA,
F_0)\stackrel{i_{*0}}{\longrightarrow} H^0(\huaA,F)
\stackrel{j_{*0}}{\longrightarrow} H^0(\huaA,\overline{F})(=0)
\\\label{LongSequenceA}
&& \longrightarrow H^1(\huaA, F_0)\stackrel{i_{*1}}{\longrightarrow}
H^1(\huaA,F) \stackrel{j_{*1}}{\longrightarrow}
H^1(\huaA,\overline{F})\longrightarrow\cdots.
\end{eqnarray}
Similarly, at any $x\in M$ the Lie algebra $L_x$ has the trivial
representation on $F_{0x}$ and an induces representation on
$\overline{F}_x$. So we have another exact sequence of complexes
$$
 0\lon C^k(L_x, F_{0x})\stackrel{i}{\longrightarrow} C^k(L_x,F_x)
\stackrel{j}{\longrightarrow} C^k(L_x,\overline{F}_x)\lon  0.
$$
And it induces a long exact sequence,
\begin{eqnarray}\nonumber
&& 0\longrightarrow H^0(L_x,
F_{0x})\stackrel{\underline{i}_{*0}}{\longrightarrow} H^0(L_x,
F_{x}) \stackrel{\underline{j}_{*0}}{\longrightarrow}
H^0(L_x,\overline{F}_x)(=0)
\\\label{LongSequenceL}
&& \longrightarrow H^1(L_x,
F_{0x})\stackrel{\underline{i}_{*1}}{\longrightarrow} H^1(L_x,
F_{x}) \stackrel{\underline{j}_{*1}}{\longrightarrow}
H^1(L_x,\overline{F}_x)\longrightarrow\cdots.
\end{eqnarray}
Moreover, there is a series of vertical arrows $\Upsilon^k$ between
(\ref{LongSequenceA}) and (\ref{LongSequenceL}), namely the
localization maps with respect to the preceding representations,
such that the diagram commute. Here we pick out a part of the
diagram as follows
\begin{equation}\label{diagram2times3}
\begin{array}{ccccccccc}
0 & \longrightarrow & H^1(\huaA,F_0)
&\stackrel{i_{*1}}{\longrightarrow} & H^1(\huaA,F) &
\stackrel{j_{*1}}{\longrightarrow} &
H^1(\huaA,\overline{F}) & \longrightarrow & \cdots \\
~&~ &\downarrow \Upsilon_{0x}^1 &~ &\downarrow \Upsilon_x^1
&~ &\downarrow \overline{\Upsilon}_x^1 &~ &~ \\
0 & \longrightarrow & H^1(L_x,F_{0x})
&\stackrel{\underline{i}_{*1}}{\longrightarrow} & H^1(L_x,F_x) &
\stackrel{\underline{j}_{*1}}{\longrightarrow} &
H^1(L_x,\overline{F}_x) & \longrightarrow & \cdots .
\end{array}
\end{equation}

By 2) of Theorem \ref{Thm:UpsilonEqual} 
and $H^0(L_x,\overline{F}_x)=0$, we know that
$\overline{\Upsilon}_x^1$ is an injection.
\begin{cor}
$$Ker({\Upsilon}_x^1)=Ker(j_{*1})\cap i_{*1}(
Ker({\Upsilon}_{0x}^1)) \cong Ker({\Upsilon}_{0x}^1).$$
\end{cor}
\pf Suppose that $\omega\in H^1(\huaA,F)$ satisfies
$\Upsilon_x^1(\omega)=0$, then
$$\underline{j}_{*1}\circ\Upsilon_x^1(\omega)=
\overline{\Upsilon}_x^1\circ j_{*1}(\omega)=0.$$ Since
$\overline{\Upsilon}_x^1$ is an injection, we have
$j_{*1}(\omega)=0$. Hence we know that
$$\omega\in Ker(j_{*1})=Im(i_{*1})\cong H^1(\huaA,F_0).$$
Since the left square in 
(\ref{diagram2times3}) commutes, we have $\omega\in Ker(j_{*1})\cap
i_{*1}( Ker({\Upsilon}_{0x}^1))$. Conversely, given any $\omega$ as
above, there naturally holds $\Upsilon_x^1(\omega)=0$.

Since $i_{*1}$ is also an injection, we have the isomorphism in the
expression. \qed

We can of course  directly regard $H^1(\huaA,F_0)
=Ker(j_{*1})\subset H^1(\huaA,F)$ and therefore
$\Upsilon_{0x}^1=\Upsilon_{x}^1|_{H^1(\huaA,F_0)}$. The above
corollary says that to know what $Ker({\Upsilon}_x^1)$ really is, it
suffices to study $Ker({\Upsilon}_{0x}^1)$.
\begin{thm}\label{Thm:KerUpsilon2}
$H^1(L_x,F_{0x})=D^1(L_x,F_{0x})$, and
\begin{equation}\label{Eqt:KerUpsilon2}
Ker({\Upsilon}_x^1)\cong Ker({\Upsilon}_{0x}^1)=\set{ [\Omega]|~~~~~
\Omega\in D^1(\huaA,F_0), \Omega|_{L_x}=0}.
\end{equation}
\end{thm}
\pf By definition of the complex $C^k(L_x,F_{0x})$, we have
$\DiffOper F_{0x}=0$ and hence
$H^1(L_x,F_{0x})=D^1(L_x,F_{0x})$.

Given any $[\Omega]\in H^1(\huaA,F_0)$, where $\Omega\in
D^1(\huaA,F_0)$, we have
$$
\Upsilon_{0x}^1([\Omega])=[\Omega|_{L_x}]=\Omega|_{L_x}.$$
So we have the expression in 
(\ref{Eqt:KerUpsilon2}). \qed

Now we have the second way expressing the kernel of the localization
map.
\begin{thm}\label{Thm:KerUpsilonSecond}
$Ker(\Upsilon_x^1)=Im[{\rho^*}(\cdot)] \cong H^1_{deR}(M,F_0)$.
\end{thm}
\pf We first point out that, in 
(\ref{Eqt:KerUpsilon2}), an $\Omega\in D^1(\huaA,F_0)$ satisfies
$\Omega|_{L_x}=0$, for some $x\in M$, if and only if
$\Omega|_{L_y}=0$ hold for all $y\in M$, i.e., or $\Omega\in
D^1(\huaA,F_0)_0$. In fact, the first conclusion of Theorem
\ref{Thm:UpsilonEqual} claims that the kernel of the localization
map does not depend on the choice of the points: $Ker(\Upsilon_x^1)=
Ker(\Upsilon_y^1)$. So the set
described in 
(\ref{Eqt:KerUpsilon2}) does not depend on the choice of $x$.

Combining Theorem 
\ref{Thm:KerUpsilon2}
with these facts and the correspondence given by Lemma 
\ref{Lem:Correspondence}, we know that each element in the kernel of
$\Upsilon_x^1$ must be the cohomology class of some
${\rho^*\theta}\in D^1(\huaA,F_0)_{0}$. \qed

We restate the conclusions of Corollary 
\ref{Cor:KerUpsilon1},
Theorem 
\ref{Thm:KerUpsilon2},
and Theorem 
\ref{Thm:KerUpsilonSecond}
 in the following theorem.
\begin{thm}\label{Thm:Summary}
Suppose that a transitive Lie algebroid
$(\huaA,[~\cdot~,~\cdot~]_\huaA,\rho)$ has a representation on a
vector bundle $F\lon M$, $\Rep: \huaA\lon \CDO(F)$. Let $p:
\widetilde{M}\lon M$ be a universal covering map. The pull back Lie
algebroid $\widetilde{\huaA}=p^!\huaA$ has an induced representation
on $\widetilde{F}=p^!F$, denoted by $\widetilde{\Rep}$. Write
$$F_0=\set{\nu\in F_y| y\in M, \Rep_u \nu=0, \forall u\in L_y}.$$
The Lie algebroid $\huaA$ has an induced representation on $F_0$,
also denoted by $\Rep$. Let $\bar{\Rep}: TM\lon \CDO(F_0)$ be the
reduced (flat) connection of $F_0$ coming from $\Rep$. Then for each
$\Omega\in D^1(\huaA,F)$, $x\in M$, the following six statements are
equivalent.
\begin{itemize}
\item[1)] $\delta_x\defbe \Omega|_{L_x}$ is a coboundary,
i.e., $\exists \tau \in {F_x}$, such that
$$
\delta_x(u)=\underline{\Rep}_u{\tau},\quad \forall u\in L_x.
$$
\item[2)] $\delta_y\defbe \Omega|_{L_y}$ is a coboundary, for
every $y\in M$.
\item[3)] The pull back cochain $p^{1*}\Omega\in D^1(\widetilde{A},\widetilde{F})$
is a coboundary, i.e., there exists some $\widetilde{\mu} \in
\Gamma(\widetilde{F})$, such that
$$
\Omega(A)=\widetilde{\Rep}_{(\widetilde{y},A)}\widetilde{\mu},\quad
\forall A\in \huaA_y.
$$
Here $\widetilde{y}\in\widetilde{M}$ satisfies $p(\widetilde{y})=y$.
\item[4)] There exists an $\Omega_0\in D^1(\huaA,F_0)$,
$\Omega_0|_{L_x}=0$, such that $\Omega$ and $\Omega_0$ are
homologic.
\item[5)] There exists an $\Omega_0\in D^1(\huaA,F_0)$,
$\Omega_0|_{L}=0$,
such that $\Omega$ and $\Omega_0$ are homologic.
\item[6)] There exist ${\theta}\in D^1(TM, F_0)$ and
$\mu\in\Gamma(F)$, such that
$$\Omega=\rho^*{\theta}+\DiffOper\mu.$$
\end{itemize}
\end{thm}

\begin{cor}\label{CorImYCong}
For any $x$, $y\in M$, the image of localizations $\Upsilon^1_x$ and
$\Upsilon^1_y$ are naturally isomorphic. That is, the isomorphism
$J$ in Diagram {\rm(\ref{Eqt:UpsilonCommute0})} naturally defines an
isomorphism $J$: $H^1(L_y,{F_y})\cong H^1(L_x,{F_x})$.
\end{cor}
\pf Consider the commute diagram which follow from Diagram
{\rm(\ref{Eqt:UpsilonCommute0})},
$$\xymatrix@R=0.5cm{
                &         Im\Upsilon_y^1 \ar[dd]^{J}     \\
H^1({\huaA},{F})\quad\quad
 \ar[ur]^{{\Upsilon_y^1}} \ar[dr]_{{\Upsilon_x^1}}                 \\
                &         Im\Upsilon_x^1 }
                $$
Since it is proved in Theorem \ref{Thm:Summary} that
$Ker(\Upsilon_x^1)=Ker(\Upsilon_y^1)$, the map $J$ in the above
diagram must be an isomorphism and naturally defined. \qed

\section{Application of
the localization theories for principal bundles and their associated
bundles}

The remaining part of this paper is devoted to apply the preceding
localization theories to that of principal bundles and their
associated vector bundles. The idea originally appeared in
   \cite{CL} and in what follows, we will recover some results in that text.

We first recall basic facts about principal bundles and the
associated bundles. Let $(P,\stackrel{\pi}{\rightarrow},M;G)$ be a
principal bundle with structure group $G$ (a Lie group) on the base
manifold $M$. We always assume that $G$ freely acts on $P$ to the
\emph{right}.

The action of $G$ on $P$ naturally lifts to an action on $TP$. We
denote the orbit of $w\in T_pP$ by $[w]$ and quotient manifold by
$\frac{TP}{G}$. Since this action is free, $\frac{TP}{G}$ admits a
vector bundle structure with base $M$, and bundle projection $q: [w]
\mapsto \pi(p)$. Sections of $\frac{TP}{G}$ can be regarded as
vector fields on $P$ which are $G$-invariant:
$$
\Gamma(\frac{TP}{G})=\set{U\in {\mathcal{X}}(P) | U_{p.g}=R_{g*}U_p,
\quad\forall p\in P,g\in G}.
$$
It follows that $\Gamma(\frac{TP}{G})$ has an induced Lie bracket
structure transferred from ${\mathcal{X}}({P})$. Besides, the
tangent map $\pi_*$ can also be transferred to $\frac{TP}{G}\lon
TM$. Thus, $(\frac{TP}{G},\pi_*,M)$ is a Lie algebroid which is
transitive (known as the \textit{gauge algebroid}).

The Artiyah sequence for this algebroid is as follows
\begin{equation}\label{SequenceAtiyahPrincipal}
 0\lon  \frac{T^\bot P}{G}\stackrel{i}{\longrightarrow} {\frac{TP}{G}}
\stackrel{\pi_*}{\longrightarrow} TM\lon  0.
\end{equation}

Here by $T^\bot P$ we denote the collection of all vertical vectors,
which are in fact the set
$$\set{\widehat{{\frkx}}_p\defbe \frac{d}{dt}|_{t=0} p.\exp t{\frkx};
\quad p\in P,{\frkx}\in \lieG}.$$

If we also have a right action of $G$ on a vector space $V$, then we
obtain the associated bundle (by $G$ acting diagonally on $P\times
V$):
$$F=\frac{P\times V}{G}.$$
Elements in $F$ are of the form $[p,v]$, for $p\in P$, $v\in V$.
Most importantly, sections of $F$ can be naturally regarded as
$V$-valued, $G$-equivariant functions on $P$. I.e.,
$$
\Gamma(F)\cong \CWGPV=\set{\mu\in C^\infty(P,V)|
\mu(p.g)=\mu(p).g,\quad\forall p\in P,g\in G}.
$$
In fact, for any $\mu\in C^\infty(P,V)$, it can be regarded as a
$V$-valued function on $P$ given by
$$
\mu(p)\defbe v, \quad \mbox{ such that }\ [p,v]=\mu(\pi(p)).
$$
And conversely, a $V$-valued, $G$-equivariant function $\mu\in
\CWGPV$ corresponds to the section of $F$ given by
$$
\mu(x)\defbe [p,\mu(p)],\quad\mbox{ by choosing an arbitary } p \in
\pi\inverse(x),\quad\forall x\in M.
$$

The Lie group $G$ has the canonical adjoint action on $\LieG$, and
hence $G$ has a right action on $\LieG$ defined by
$$
{\frkx}.g\defbe Ad_{g\inverse}{\frkx},\quad\forall {\frkx}\in \LieG,
g\in G.
$$

The adjoint Lie algebra bundle $\frac{T^\bot P}{G}$ of the
gauge Lie algebroid in Sequence 
(\ref{SequenceAtiyahPrincipal}) is indeed the associated bundle
$$
\frac{P\times Lie(G)}{G}.
$$
In fact, given any $\LieG$-valued $G$-equivariant function $\kappa$:
$P\lon \LieG$, which satisfies
$$\kappa(p.g)=\kappa(p).g=Ad_{g\inverse}\kappa(p),\quad
\forall p\in P, g\in G,$$ it corresponds to a vertical vector field
on $P$:
\begin{equation}\label{widehatkappa}
\widehat{\kappa}|_{p}\defbe
\widehat{\kappa(p)}_p=\frac{d}{dt}|_{t=0} p.\exp t{\kappa(p)},\quad
\forall p\in P.
\end{equation}
One can directly check that $\widehat{\kappa}$ is a $G$-invariant
vector field  on $P$.

The ring of smooth functions $\CWM$ can also be regarded as
$G$-invariant functions on $P$:
$$
\CWM\cong \pi^*\CWM=\set{f\in C^\infty(P)| f(p.g)=f(p), \quad\forall
p\in P,g\in G}.
$$

There is a standard representation of the gauge algebroid
$\frac{TP}{G}$ on $F$, defined by
$$\Rep_{U}(\mu)\defbe U(\mu),\quad\forall U\in\Gamma(\frac{TP}{G}),
\mu\in \CWGPV.$$

Let $\Omega^{1,G}(P,V)$ denote the $V$-valued, $G$- equivariant
1-forms on $P$:
\begin{eqnarray*}
&&\Omega^{1,G}(P,V)\\
&=&\set{\omega: TP\lon V;\omega|_{p}(w).g
=\omega|_{p.g}(R_{g*}w),\forall p\in P,g\in G,
w\in T_pP}\\
&=&\set{\omega: TP\lon V; \omega(U)\in \CWGPV,\quad \forall
U\in\Gamma(\frac{TP}{G})}.
\end{eqnarray*}
\begin{lem}
There is a canonical pull back morphism of $\CWM$-modules
$$\pi^*: \Gamma(Hom(TM,F))\lon \Omega^{1,G}(P,V),
\quad \vartheta\mapsto \pi^*(\vartheta),$$ where
$$\pi^*(\vartheta): \quad w\mapsto \vartheta(\pi_*w)(p), \quad \forall w\in T_pP.$$
Moreover, $\pi^*$ is injective.
\end{lem}
\pf To see that $\pi^*(\vartheta)$ is a $G$- equivariant 1-form, we
suppose that for $w\in T_pP$, $R_{g*}w\in T_{p.g}P$, $x=\pi(p)\in M$
and
$$X=\pi_*(w)=\pi_*(R_{g*}w)\in T_xM,$$
one has $\vartheta(X)=[p,v]=[p.g,v.g]$. Then,
\begin{eqnarray*}
&&\pi^*(\vartheta)|_p(w).g=\vartheta(X)(p).g\\
&=&v.g=\vartheta(X)(p.g)=\pi^*(\vartheta)|_{p.g}(R_{g*}w).
\end{eqnarray*}
It is easy to see that $\pi^*(\vartheta)$ being a zero $1$-form
implies that $\vartheta$ is zero. \qed

We also introduce the notation $\Omega^{k,G}(P,V)$ ($k\geqslant 1$)
to denote the $V$-valued $G$-equivariant $k$-forms on $P$. I.e.,
\begin{eqnarray*}
&&\Omega^{k,G}(P,V)\\
&=&\set{\omega: \wedge^k TP\lon V;\omega|_{p}(w).g
=\omega|_{p.g}(R_{g*}w),\forall p\in P,g\in G,
w\in \wedge^k T_pP}\\
&=&\set{\omega: \wedge^k TP\lon V;\omega(U)\in \CWGPV,\quad \forall
U\in\Gamma(\wedge^k \frac{TP}{G})}.
\end{eqnarray*}
These $\Omega^{k,G}(P,V)$ together with $\Omega^{0,G}(P,V)=\CWGPV$
become a complex (over the ring $\CWM$), equipped with the usual
exterior differential operator $d$. It is in fact isomorphic (as
complexes) to $C^k(\frac{TP}{G},F)$. And in turn, we have
$$
H^{k,G}_{deR}(P,V)\cong H^k(\frac{TP}{G},F).
$$

\begin{defi}
Given an arbitrary $p\in P$, the localization map for the gauge
algebroid $\frac{TP}{G}$ and its associated vector bundle
$F=\frac{P\times V}{G}$ is a group morphism
$$
\Upsilon^k_p: H^{k,G}_{deR}(P,V)\lon H^k(\lieG,V).
$$
It sends the cohomology class of $\omega\in \Omega^{k,G}(P,V)$ to
the cohomology class of $\widehat{\omega}|_{p}\in D^k(G,V)$ which is
defined by
$$
\widehat{\omega}|_p ({\frkx})\defbe \omega|_p(\widehat{{\frkx}}_p),
\quad\forall {\frkx}\in Lie(G).
$$

\end{defi}
Let $G_e$ denote the subgroup of $G$ which is the connected
component of $G$ containing the unit element $e$.
Consider a sub vector space
$$V_0=\set{v\in V| v.h=v,\quad
\forall h\in G_e }.$$ One is easy to prove that $V_0$ is
$G$-invariant. Hence $G$ also has a right action on $V_0$.

\begin{lem}
For the representation of the gauge Lie algebroid
 $\frac{TP}{G}$ on $F=\frac{P\times V}{G}$,
$$F_0=\frac{P\times V_0}{G}=H^0(\frac{T^\bot P}{G},F).$$
\end{lem}
\pf Let $\kappa$: $P\lon \LieG$ be a $\LieG$-valued $G$-equivariant
function which corresponds to a $G$-invariant vertical vector field
$\widehat{\kappa}$ on $P$ given by 
(\ref{widehatkappa}). Let $\mu\in \CWGPV$ be a $V$-valued,
$G$-equivariant function which can also be regarded as an element of
$\Gamma(F)$. At any $p\in P$,
 $s\in \Real$, let $g=\exp s\kappa(p)$. Then we have
\begin{eqnarray*}
&&\widehat{\lambda}|_{p.g}(\mu)\\
&=& \frac{d}{dt}|_{t=0} \mu(p.g.\exp t\kappa(p.g))
 =
\frac{d}{dt}|_{t=0} \mu(p).g.\exp tAd_{g\inverse}\kappa(p)\\
& =&\frac{d}{dt}|_{t=0} \mu(p).\exp t\kappa(p).g=
\frac{d}{dt}|_{t=0} \mu(p).\exp t\kappa(p).g\\
& =&\frac{d}{dt}|_{t=0} \mu(p).\exp (t+s)\kappa(p)=
\frac{d}{dt}|_{t=s} \mu(p).\exp t\kappa(p).
\end{eqnarray*}
Hence we know that $\mu\in \Gamma(F_0)$ implies
$$
\frac{d}{dt} \mu(p)\exp t{\kappa(p)}=0,\quad \mbox{i.e.,} \ \
\mu(p)\exp t{\kappa(p)}\equiv \mu(p),\ \forall t\in \Real.
$$
Since $\kappa$ is arbitrary and $G_e$ is generated by elements of
the form $\exp t{\frkx}$, ${\frkx}\in \LieG$, $\mu(p)$ is an element
of $V_0$. This shows that the function $\mu$ takes values in $V_0$.
Conversely, if $\mu\in \Gamma(F_0)$, then it obviously satisfies
$\widehat{\kappa}(\mu)=0$, $\forall \widehat{\kappa}\in \frac{T^\bot
P}{G}$. \qed

Using Theorem 
\ref{Thm:KerUpsilonSecond}, we have the following conclusion
describing the kernel of the localization map $\Upsilon^1_p$.
\begin{thm}
$$Ker(\Upsilon^1_p)=\pi^*H^1_{deR}(M,F_0)\cong H^1_{deR}(M,F_0)\,.$$
Here $\pi^*: H^1_{deR}(M,F_0)\lon H^{k,G}_{deR}(P,V)$ is given by
$$
[\theta]\mapsto [\pi^*(\theta)],\quad\forall \theta\in D^1(TM, F_0),
$$
and it is an injection.
\end{thm}

And one may restate the above theorem into the following
form, analogue to that of Theorem 
\ref{Thm:Summary}.
\begin{thm}With the preceding notations, let
$w\in \Omega^{1,G}(P,V)$ be a close $1$-form. For any $p\in P$, we
have five equivalent statements:
\begin{itemize}
\item[1)] $\omega|_{T^\bot_p P}$ is a coboundary, i.e.,
$\exists v\in V$ such that
$$\omega(\widehat{{\frkx}}_p)=\frac{d}{dt}|_{t=0} v.\exp t{\frkx},\quad \forall {\frkx}\in\lieG.$$

\item[2)] $\omega|_{T^\bot_q P}$ is a coboundary, for all $q\in P$.
\item[3)] There exists a closed $1$-form $\omega_0\in \Omega^{1,G}(P,V_0)$,
$\omega_0|_{T^\bot_p P}=0$, such that $\omega_0$ and $\omega$ are
homologic, i.e.,
$$\omega=
\omega_0+d \mu,\quad\mbox{for some }\mu\in\CWGPV.
$$
\item[4)] There exists a closed $1$-form $\omega_0\in \Omega^{1,G}(P,V_0)$,
$\omega_0|_{T^\bot P}=0$, such that $\omega_0$ and $\omega$ are
homologic.
\item[5)] For some closed $1$-form $\theta\in D^1(TM, F_0)$
and $\mu\in \CWGPV$, $$\omega=\pi^*\theta+ d\mu.$$
\end{itemize}
\end{thm}

In particular, we conclude that:
\begin{cor} If $M$ is simply connected (or $V_0=0$), then
the localization map $\Upsilon^1_p$ is an injection. In other words,
any one of the five statements in the above theorem implies that
$\omega$ is a coboundary, i.e., $\omega=d\mu$ for some
$\mu\in\CWGPV$.
\end{cor}

\section{Transitive Lie Bialgebroids}
In this section ,   we recall some results on the  structure of
transitive Lie bialgebroids in  \cite{CL} as an another application
of our localization theory.

\noindent\textbf{$\bullet$ Lie bialgebroids.} A Lie bialgebroid is a
pair of Lie algebroids (${\cal A}$, ${\cal A}^*$) satisfying the
following compatibility condition
\begin{equation}
\label{eq:1} d_{*}[A, B]_\huaA=[d_{*}A, B]_\huaA+[A, d_{*}B]_\huaA,
\ \ \forall A , B \in \Gamma ({\cal A}),
\end{equation}
where the differential $d_*$ on $\Gamma(\wedge^{\bullet}{\cal A})$
comes from the Lie algebroid structure on ${\cal A}^*$ (see
\cite{K-S:1994}, \cite{MackenzieX:1994} for more details). Of
course, one can also denote a Lie bialgebroid by the pair $({\cal
A},\ d_*)$, since the
 anchor $\rho_* : {\cal A^*} \lon TM $ and
the Lie bracket $[~\cdot~, ~\cdot~]_* $ on the dual bundle are
defined by $d_*$ as follows: $\rho^*_* (df) = d_*f, ~\forall f \in
\cm$ and for all $ A \in \ga
    ,~ \xi , \eta \in \Gamma ({\cal A}^{*})$,
$$\langle[\xi , \eta ]_*, A \rangle= \rho_*(\xi)\langle\eta, A\rangle - \rho_*(\eta)\langle\xi,
A \rangle    -  d_*A (\xi, \eta).$$ 

Again we suppose that the Lie algebroid $\huaA$ is transitive.
Recall the adjoint representation of $\huaA$ on $L$ defined in
(\ref{adjointRep}). In this paper, we also consider the adjoint
representation of $\huaA$ on $L\wedge L$ associated to that on $L$,
and we write $L^2$ for $L\wedge L$. In this case, $\Omega\in
C^1(\huaA, {L^2})$, i.e., a bundle map from $\huaA$ to $L^2$, is a
\emph{Lie algebroid 1-cocycle} if and only if
\begin{equation}\label{EqtCocycleOmegaL2}
\Omega[A,B]_\huaA
=[\Omega(A),B]_\huaA+[A,\Omega(B)]_\huaA,\quad\forall A,B\in\Gamma
(\huaA).
\end{equation}
$\Omega$ is a coboundary if $\Omega=[\mu,~\cdot~ ]_\huaA$ for some
$\mu\in\Gamma(L^2)$.

The structure of transitive Lie bialgebroids is studied in
\cite{CL}. We quote directly some of the conclusions in that text.
\begin{defi}
For a transitive Lie algebroid
$(\huaA,[~\cdot~,~\cdot~]_\huaA,\rho)$, given
$\Lambda\in\Gamma(\wedge^2 \huaA)$ and a bundle map $\Omega$:
$\huaA\lon  {L^2}$,
the pair $(\Lambda,\Omega)$ is called $\huaA$-compatible 
if $\Omega$ is a 1-cocycle and satisfies 
\begin{equation}\label{EqtDDzero}
[\frac{1}{2}[\Lambda,\Lambda]_\huaA+\Omega(\Lambda)
,~\cdot~]_\huaA+\Omega^2=0, \quad\mbox{as a map } \Gamma(\huaA)\lon
\Gamma(\wedge^3 \huaA).
\end{equation}
\end{defi}
Here $\Omega(\Lambda)$ and $\Omega^2$ make sense by means of the
extension of $\Omega$ as a derivation of  the graded bundle,
$\Omega$: $\wedge^{k} \huaA\lon \wedge^{k+1}\huaA$, $k\geq 0$. For
$k=0$, it is zero. For $k\geq 1$, it is  defined by
\begin{equation}\label{ruleofextend}
\Omega(A_1\wedge\cdots\wedge
A_k)=\sum_{i=1}^{k}(-1)^{i+1}A_1\wedge\cdots\wedge\Omega
(A_i)\wedge\cdots\wedge A_k,\ \
\end{equation}
for all $A_1\wedge\cdots\wedge A_k\in \Gamma(\wedge^k \huaA)$.

It is easy to see that if $(\Lambda,\Omega)$ is $\huaA$-compatible,
then so is the pair $(\Lambda+\nu,\Omega-[\nu,~\cdot~]_\huaA)$, for
any $\nu\in\Gamma({L^2})$. Thus, two $\huaA$-compatible pairs
$(\Lambda,\Omega)$ and $(\Lambda',\Omega')$ are called equivalent,
written $(\Lambda,\Omega)\sim(\Lambda',\Omega')$, if
$\exists\nu\in\Gamma({L^2})$, such that $\Lambda'=\Lambda+\nu$ and
$\Omega'=\Omega-[\nu,~\cdot~]_\huaA$.

\begin{thm}\label{ThmFirst}
Let $(\huaA,[~\cdot~,~\cdot~]_\huaA,\rho)$ be a {\em transitive} Lie
algebroid over $M$. Then there is a one-to-one correspondence
between Lie bialgebroids $(\huaA,d_*)$ and equivalence classes of
$\huaA$-compatible pairs $(\Lambda,\Omega)$ such that
\begin{equation}\label{EqtRelationDStarAndPair}
d_*=[\Lambda,~\cdot~]_\huaA+\Omega.
\end{equation}
\end{thm}

For a Lie bialgebra $(\f{g},\f{g}^*)$, it is obvious that one can
take $\Lambda = 0$ and $-\Omega$  as  the cobracket of $\f{g}$.
Another special case is the following.

\begin{cor}If $\Omega$:
$\huaA\lon  {L^2}$ is a 1-cocycle and satisfies $\Omega^2=0$, as a
map $\huaA\lon \wedge^3 \huaA$, then $(\huaA,\Omega)$ is a Lie
bialgebroid. In this case, the anchor $\rho_*$ of $\huaA^*$ is zero
and $\huaA^*$ is a bundle of Lie algebras whose bracket is defined
by
$$\langle[\xi , \eta ]_*, A
\rangle=-\langle\Omega(A),\xi\wedge\eta\rangle,$$ for all
$\xi,\eta\in \huaA^*$ and $A\in\huaA$.
\end{cor}

\begin{cor}\label{CorPullBackLieBialgebroid}
Let $(\huaA,d_*)$ be a transitive Lie algebroid over $M$ and suppose
that $d_*=[\Lambda,~\cdot~]_\huaA+\Omega$ is given as in
{\rm(\ref{EqtRelationDStarAndPair})}. Let $p: \widetilde{M}\lon M$
be a covering. Then the pull back bundle
$\widetilde{\huaA}=p^!\huaA$ is also a transitive Lie algebroid. For
the pull back section
$\widetilde{\Lambda}\in\Gamma(\widetilde{\huaA})$ and the pull back
bundle map
$$
\widetilde{\Omega}: \quad \widetilde{\huaA}\lon \widetilde{L}^2,
$$
let
$$
\widetilde{d}_*=[\widetilde{\Lambda},~\cdot~]_{\widetilde{\huaA}}+\widetilde{\Omega}.
$$
Then $(\widetilde{\huaA},\widetilde{d}_*)$ is also a Lie bialgebroid
over $\widetilde{M}$.
\end{cor}

It is known that, for any section $\Lambda\in
\Gamma(\wedge^2\huaA)$, 
one can define a bracket on $\Gamma ({\cal A}^* )$  by
$$[\xi , \eta ]_{\Lambda}=L_{\rr \xi}\eta -L_{\rr \eta}\xi -d<\rr
\xi , \eta >.$$ With the bracket defined above and anchor map
$$\rho_* \defbe \rho\circ \rr :\quad {\cal A^*}  \lon TM,$$
the dual bundle ${\cal A}^*$ becomes  a Lie algebroid if and only if
$[X, [\rt, \rt]_{\huaA}]_{\huaA} = 0$, \,$\forall X \in \Gamma
({\cal A})$ (\cite{ExactLieBlg}, Theorem 2.1). In this situation,
the induced differential on $\Gamma(\wedge^{\bullet}{\cal A})$ has
the form, $d_* = [\rt, ~\cdot~]_\huaA$, and clearly satisfies
compatibility condition (\ref{eq:1}). The Lie bialgebroid arising in
this way is called a {\em coboundary (or exact) Lie bialgebroid}
 \cite{ExactLieBlg}. In the particular case where $[\rt, \rt]_{\huaA}
= 0$, the Lie bialgebroid is called {\em triangular} \cite{MX2}.

 By our
definition of $\huaA$-compatible pairs, the pair corresponding to a
coboundary Lie bialgebroid can be chosen to be $(\Lambda,0)$ or,
equivalently,  $(\huaA,\huaA^*)$ is a coboundary Lie bialgebroid if
and only if the second element of the $\huaA$-compatible pair
$\Omega\in C^1(\huaA,{L^2})$ is a coboundary. Therefore, to deal
with coboundary Lie bialgebroids, one first needs to study the
properties of  Lie algebroid 1-cocycles.

\begin{cor} With the same assumptions as in Theorem \ref{ThmFirst},
if $Rank(\huaA)=1$, then $d_*=0$. That is, $\huaA^*$ admits trivial
Lie algebroid structures.
\end{cor}

\begin{cor} With the same assumptions as in Theorem \ref{ThmFirst},
if $Rank(L)=1$, then $d_*=[\Lambda,~\cdot~]$. That is, $(\huaA,d_*)$
is coboundary.
\end{cor}

The following theorem follows directly from 2) of Theorem
\ref{Thm:UpsilonEqual}.
\begin{thm}\label{ThmOmegaDelta}
Suppose that a transitive Lie algebroid $\huaA$ satisfies  {\em
\textbf{one}}  of the following conditions
\begin{itemize}
\item[\rm 1)] $H^0(\f{g},\f{g}^2)=0$, where $\f{g}= L_x$, for some $x \in
M$;
\item[\rm 2)] $M$ is simply connected.
\end{itemize}
Then $\Omega\in C^1(\huaA,L^2)$ is coboundary if and only if
$\delta_x\defbe \Omega|_{\f{g}}$ is coboundary.

In particular, if $H^1(\f{g},\f{g}^2)=0$, any Lie bialgebroid
$(\huaA,\huaA^*)$ is coboundary.
\end{thm}

It is a well known result that for any nontrivial representation of
a semi-simple Lie algebra $\f{g}$ on some vector space $V$, the
cohomology groups $H^0(\f{g},V)$ and $H^1(\f{g},V)$ are both zero.
So we conclude:
\begin{cor}Let $\huaA$ be a transitive Lie algebroid and let
$\f{g}= L_x$ be the fiber type of $L$. If $\f{g}$ is semi-simple and
its adjoint representation on $\f{g}^2$ is not trivial, then any Lie
bialgebroid $(\huaA,\huaA^*)$ is coboundary.
\end{cor}

We also have the following corollaries which follow from Theorem
\ref{Thm:Summary}. Note that now
$$
L^2_0=\set{\nu\in L^2_y| y\in M, [u, \nu]_L=0, \forall u\in L_y}
$$
is a vector bundle over $M$ which has a natural flat connection.
\begin{cor} Suppose that a transitive Lie bialgebroid $(\huaA,d_*)$ satisfies
$H^1(\f{g},\f{g}^2)=0$, where $\f{g}= L_x$ for some $x \in M$. Then
\begin{itemize}
\item[1)]
for a universal covering $p$: $\widetilde{M}\lon M$, the pull back
Lie bialgebroid $(\widetilde{\huaA},\widetilde{d}_*)$ given by
Corollary {\rm\ref{CorPullBackLieBialgebroid}} is coboundary;
\item[2)]
the compatible pair corresponding to $(\huaA,d_*)$ can be chosen to
be $(\Lambda,{\rho^*\theta})$, for some closed $L^2_0$-coefficient
$1$-form ${\theta}\in D^1(TM, L^2_0)$, where
$$
{\rho^*\theta}: \quad A\mapsto \theta(\rho(A)),\quad\forall
A\in\huaA.
$$
\end{itemize}
\end{cor}
Note that in this case, $\Omega|_{L}={\rho^*\theta}|_{L}$ is
trivial, and hence $\Omega^2=0$. So the compatible relation given in
Equation (\ref{EqtDDzero}) becomes
$$
[\frac{1}{2}[\Lambda,\Lambda]_\huaA+({\rho^*\theta})\lrcorner\Lambda
,~\cdot~]_\huaA=0, \quad\mbox{as a map } \Gamma(\huaA)\lon
\Gamma(\wedge^3 \huaA).
$$

\end{document}